
\input amstex.tex
\documentstyle{amsppt}
\def\DJ{\leavevmode\setbox0=\hbox{D}\kern0pt\rlap
{\kern.04em\raise.188\ht0\hbox{-}}D}
\footline={\hss{\vbox to 2cm{\vfil\hbox{\rm\folio}}}\hss}
\nopagenumbers
\font\ff=cmr8
\def\txt#1{{\textstyle{#1}}}
\baselineskip=13pt
\def\hf{{\textstyle{1\over2}}}
\def\a{\alpha}\def\b{\beta}

\def\d{{\,\roman d}}
\def\e{\varepsilon}
\def\f{\varphi}

\def\G{\Gamma}
\def\k{\kappa}
\def\s{\sigma}
\def\t{\theta}
\def\ti{{\tilde a}(n)}
\def\Ti{{\tilde a}*{\tilde a}(n)}
\def\={\;=\;}

\def\D{\Delta}
\def\no{\noindent}
\def\R{\Re{\roman e}\,} 
\def\z{\zeta}
\def\e{\varepsilon}
\def\D{\Delta}
\def\no{\noindent}

\def\st{\sigma_k^*}
\font\teneufm=eufm10
\font\seveneufm=eufm7
\font\fiveeufm=eufm5
\newfam\eufmfam
\textfont\eufmfam=\teneufm
\scriptfont\eufmfam=\seveneufm
\scriptscriptfont\eufmfam=\fiveeufm
\def\mathfrak#1{{\fam\eufmfam\relax#1}}

\font\tenmsb=msbm10
\font\sevenmsb=msbm7
\font\fivemsb=msbm5
\newfam\msbfam
\textfont\msbfam=\tenmsb
\scriptfont\msbfam=\sevenmsb
\scriptscriptfont\msbfam=\fivemsb
\def\Bbb#1{{\fam\msbfam #1}}

\def \NN {\Bbb N}
\def \CC {\Bbb C}

\def \RR {\Bbb R}
\def \ZZ {\Bbb Z}

\def\rightheadline{{\hfil{\ff
Mean values of zeta-functions}\hfil\tenrm\folio}}

\def\leftheadline{{\tenrm\folio\hfil{\ff
Aleksandar Ivi\'c }\hfil}}
\def\emptyheadline{\hfil}
\headline{\ifnum\pageno=1 \emptyheadline\else
\ifodd\pageno \rightheadline \else \leftheadline\fi\fi}

\topmatter
\title ON MEAN VALUES OF SOME ZETA-FUNCTIONS
IN THE CRITICAL STRIP \endtitle

\author   Aleksandar Ivi\'c
\bigskip
{\sevenbf Dedicated to the memory of Robert Rankin}
\endauthor
\address{
Aleksandar Ivi\'c, Katedra Matematike RGF-a
Universiteta u Beogradu, \DJ u\v sina 7, 11000 Beograd,
Serbia (Yugoslavia).}
\endaddress
\email {\tt aivic\@matf.bg.ac.yu,
aivic\@rgf.bg.ac.yu} \endemail

\abstract
{For a fixed integer $k \ge 3$, and fixed $\hf < \s < 1$ we consider
$$
\int_1^T|\z(\s+it)|^{2k}\d t = \sum_{n=1}^\infty d_k^2(n)n^{-2\s}T
+ R(k,\s;T),
$$
where $R(k,\s;T) = o(T)\;(T\to\infty)\,$ is the error term in the above
asymptotic formula. Hitherto the sharpest bounds for $R(k,\s;T)$ are derived
in the range $\min(\b_k,\st) < \s < 1$. We also obtain new mean
value results for the zeta-function of holomorphic cusp forms
and the Rankin-Selberg series.}
\endabstract

\endtopmatter

\heading{\bf1. Introduction}
\endheading

The  aim of this paper is to provide asymptotic formulas for the $2k$--th
moment of the Riemann zeta-function $\z(s)$  and some related Dirichlet
series  in the so-called ``critical strip" $\hf < \s = \R s < 1$.
For the zeta-function our results are relevant when $k\ge3$ is a
fixed integer, where henceforth  $s = \s + it$ will denote a complex variable.
Mean  values of $\z(s)$ on the ``critical line" $\s = \hf$ behave differently
(see e.g., [4]), while the problem of mean values for $0 < \s < \hf$
can be reduced to the range $\hf < \s < 1$ by means of the functional
equation for $\z(s)$, namely
$$
\z(s) = \chi (s)\z(1-s),\quad \chi(s) = 2^s\pi^{s-1}\sin(\hf \pi s)\G(1-s).
$$
Mean values of $\z(s)$ for $\hf < \s < 1$ in the cases $k = 1$
and $k = 2$ have been extensively studied, and represent one of the
central themes in zeta-function theory. One has (see [10, Theorem 2])
$$\eqalign{
\int_1^T|\z(\s+it)|^2\d t &= \z(2\s)T + {\z(2\s-1)\G(2\s-1)\over1-\s}
\sin(\pi\s)T^{2-2\s} \cr&+ O(T^{2(1-\sigma)/3}\log^{2/9}T)
\qquad(\hf < \sigma \le 1),\cr}\eqno(1.1)
$$
and (see [7, Theorem 2])
$$
\int_1^T|\z(\s+it)|^4\d t = {\z^2(2\s)\over\z(4\s)}T + O(T^{2-2\s}\log^3T)
\qquad(\hf < \sigma \le 1),\eqno(1.2)
$$
which are the sharpest hitherto published asymptotic formulas valid in the
whole range $\hf < \sigma \le 1$. These results have been obtained
by special methods, and cannot be generalized to higher moments.
The formula for the general $2k$--th moment of $\z(s)$ can be
conveniently written (cf. [4, Chapter 8]) as
$$
\int_1^T|\z(\s+it)|^{2k}\d t = \sum_{n=1}^\infty d_k^2(n)n^{-2\s}T
+ R(k,\s;T), \quad R(k,\s;T) = o(T),\eqno(1.3)
$$
where $\hf < \s_0(k) \le \s \le 1,\,T\to\infty$, and
the arithmetic function $d_k(n)$ denotes, as usual,
the number of ways $n$ may be written
as a product of $k$ factors (so that $d_k(n)$ is generated by $\z^k(s)$,
and $d_2(n) = d(n)$ is the number of divisors of $n$).
In [4, Chapter 8] it was proved that
$$
R(k,\s;T) \;\ll_\e\;T^{{2-\s-\st\over2-2\st}+\e}\qquad(\st < \s < 1),
\eqno(1.4)
$$
where henceforth $\e$ denotes arbitrarily small constants,
not necessarily the same ones at each occurrence, and
$\st$  is the infimum of $\s^*\,(\ge \hf)\,$ for which
$$
\int_1^T|\z(\s^* + it)|^{2k}\d t \;\ll_\e\; T^{1+\e}
$$
holds for any given $\e$. Writing further the bounds for $R(k,\s;T)$ as
$$
R(k,\s;T) \;\ll_\e\; T^{c_k(\s)+\e}
$$
and using the known bounds for $\st$ when $3\le k \le 6$, it follows from
(1.4) that we have
$$\eqalign{&
c_3(\s) = {17-12\s\over10} \quad({\txt{7\over12}} < \s < 1),\quad
c_4(\s) = {11-8\s\over6} \quad({\txt{5\over8}} < \s < 1),\cr&
c_5(\s) = {79-60\s\over38} \quad({\txt{41\over60}} < \s < 1),\quad
c_6(\s) = {9-7\s\over4} \quad({\txt{5\over7}} < \s < 1).\cr}\eqno(1.5)
$$
As indicated in [4], explicit values for $c_k(\s)$ could be given
for any fixed $k > 1$, but the expressions in general would be cumbersome,
so only explicit values were given for $2 \le k \le 6$. The point of
(1.3)-(1.5) lies in the  fact that each value of $c_k(\s)$ satisfies
$c_k(\s) < 1$ (i.e., when (1.3) becomes a true asymptotic formula),
precisely for  the range given in (1.5).
However, as $\s$ approaches 1, the values of
$c_k(\s)$ become rather poor and they do not tend to zero, as one expects.

\medskip
The problem of mean values of a Dirichlet series $F(s)$ (in this context
$2k$--th moments of $F(s)$ can be regarded simply as the mean square
of $F^k(s)\;(k\in\NN)$) can be treated in various degrees of generality.
Here we shall mention only the classes of   Dirichlet series
treated by Chandrasekharan--Narasimhan
(see [2], [3]), Perelli [14], Richert [17] and Selberg [18]. Recently
S. Kanemitsu et al. obtained in [12] a mean value theorem for a general
class of Dirichlet series possessing a functional equation with
multiple gamma-factors. The merit of their result, which is in part
based on ideas of Matsumoto [13], is a relatively
good value of the exponent in the error term as $\s$ approaches the
abscissa of absolute convergence of the Dirichlet series in question.
In particular, the result of [12] can be applied to higher power moments
of $\z(s)$. In this case in the notation of [12] one has
$$
\alpha = 0,\;\mu = \nu = 1,\; \alpha_1 = 0,\;\gamma_1 = \hf, \;\beta_1 = \hf,
\; H = 1,\;\eta = \hf.
$$
Their Theorem 4 gives then, in the notation of (1.3),
$$
R(k,\s;T) \;\ll_\e\;T^{{3k(1-\sigma)\over k+2-k\sigma}+\e}\eqno(1.6)
$$
for $k \ge 2$ and
$$
1 - {1\over k} + \varepsilon \le \sigma \le 1.\eqno(1.7)
$$
When (1.6) is compared with (1.3)--(1.4) it transpires that it holds
for a poorer range, but the exponent in the error term is much sharper
as $\s$ grows, and it tends to 0 as $\s\to1-0$, as one expects.

In what follows we may assume $\s < 1$, since we have the
asymptotic formula
$$
\int_1^T|\z(1+it)|^{2k}\d t = \sum_{n=1}^\infty |d_k(n)|^2n^{-2}T
+ O((\log T)^{|k|^2}),\eqno(1.8)
$$
which was proved in [1]. In (1.8) one can take $k\in\CC$  arbitrary,
but fixed. Thus (1.8), obtained by a special method that cannot
be adapted to the range $\s < 1$,
yields a better error term than the one obtainable from any of
the previous bounds (1.1)--(1.7).

The plan  of the paper is as follows. In Section 2 we shall formulate
the results (Theorem 1 and Theorem 2) concerning the higher moments
of $\z(s)$, the proofs of which will be given in Section 3. In
Section 4 we shall deal with the mean value of the Rankin-Selberg series,
and in Section 5 with the mean values of the zeta-function of holomorphic
modular forms and its square.

\heading{\bf2. Higher moments of the zeta-function}
\endheading

The aim of this section is to furnish new bounds for $R(k,\s;T)$,
which will improve both (1.4) and (1.6). We shall
formulate now our results, with the remark that Theorem 2 is based on the
use of the defining property of $\st$ and it gives good bound for
$R(k,\sigma;T)$  when $\s$ is close to $\st$. Theorem 1, on the other
hand, is derived by using the values of the constant $\b_k$ in the mean
square estimates for the divisor problem. Namely we let, as usual,
$$
\b_k = \inf\left\{\,b_k \;(\ge 0)\;:\;
\int_1^x\D^2_k(y)\d y \;\ll\; x^{1+2b_k}\right\},\eqno(2.1)
$$
where $\D_k(x)$ is the error term in the asymptotic formula for the
summatory function of $d_k(n)$ (cf. (3.1)).
Theorem 1 will provide good results for values of $\s$ close to 1.
Results of similar type for the general case and the case of the
Rankin-Selberg series can be found in [12] and [13]. However in the
proof of Theorem 1 we shall avoid
using the Cauchy-Schwarz inequality and therefore obtain a sharper
value of the exponent than we would obtain by following the ideas
of [12] and [13].

\medskip
THEOREM 1. {\it For fixed $\s$ satisfying $\max(\b_k,\hf) < \s < 1$ and every
fixed integer $k \ge 3$, we have}
$$
R(k,\s;T) \;\ll_\e\;T^{{2(1-\sigma)\over 1-\b_k}+\e}.\eqno(2.2)
$$

\medskip
THEOREM 2. {\it For fixed $\s$ satisfying $\st < \s < 1$ and every
fixed integer $k \ge 3$, we have}
$$
R(k,\s;T) \;\ll_\e\;T^{{2(1-\sigma)\over 2-\st-\s}+\e}.\eqno(2.3)
$$

\medskip{\bf Remark 1.}
Note that (2.2) improves (1.6).
Namely we have
$$
{2\over1-\b_k}\;\le\;{3k\over k+2-k\s}
$$
for
$$
2k\s \;\ge\;4-k+3k\b_k.
$$
But from (1.7) it follows that
$$
2k\s > 2k -2 \ge 4 - k + 2k\b_k
$$
for
$$
\b_k \;\le 1 - {2\over k}\qquad(k \ge 3). \eqno(2.4)
$$
Equality in (2.4) holds only for $k = 3$, since $\b_3 = {1\over 3}$ (see [4]).
But we have $\b_4 = {3\over8}$ and $\b_k \le (k-1)/(k+2)$ for $k \ge 4$
(see [17]), hence in (2.4) we have strict inequality for $k > 3$.
This means that (2.2) improves both the exponent
of the error term in (1.6), and at the same time it holds in a wider interval
than the one given by (1.7).

We also note that
$$
{2-2\s\over2-\st-\s} \;\le\; {2-\s - \st\over2 - 2\st}
$$
is equivalent to
$$
2\s\st \le \s^2 + (\st)^2,
$$
which is obvious.  This means that  (2.3) of Theorem 2
improves (1.4) in the whole range $\st < \s < 1$.

\medskip

\heading{\bf3. Proof of Theorem 1 and Theorem 2}
\endheading

\no
We write as usual, for $k\in\NN$,
$$
D_k(x) := \sum_{n\le x}d_k(n) = xP_{k-1}(\log x) + \D_k(x),\eqno(3.1)
$$
where $P_{k-1}(y)$ is a polynomial of degree $k-1$, whose coefficients (which
depend on $k$) may be explicitly evaluated.
Using the Stieltjes integral representation and (3.1) we have,
for $1 \ll X \ll T^C\,(C>0),\,\s > 1$ and $k\ge 2$ a fixed integer,
$$\eqalign{
\z^k(s) &= \sum_{n\le X}d_k(n)n^{-s} + \int_X^\infty x^{-s}\d D_k(x)\cr&
= \sum_{n\le X}d_k(n)n^{-s} + \int_X^\infty x^{-s}
(Q_{k-1}(\log x)\d x + \d\D_k(x)),\cr&}\eqno(3.2)
$$
where $Q_{k-1} = P_{k-1} + P'_{k-1}$. From the definition (2.1) of $\b_k$
it follows that, for any given $Y \gg 1$,
there exists $X \in [Y,\,2Y]\,$ such that
$$
\D_k(X) \ll_\e X^{\b_k+\e}.\eqno(3.3)
$$
Henceforth we assume that $X$ is chosen in such a way that it satisfies,
besides $1 \ll X \ll T^C$,
also the bound in (3.3). Repeated integration by parts yields
$$\eqalign{&
\int_X^\infty x^{-s}Q_{k-1}(\log x)\d x =
{X^{1-s}\over s-1}Q_{k-1}(\log X) + {1\over s-1}
\int_X^\infty x^{-s}Q'_{k-1}(\log x)\d x\cr&
= \ldots = X^{1-s}\left({Q_{k-1}(\log X)\over s-1}
+ {Q'_{k-1}(\log X)\over (s-1)^2} + \ldots +
{Q^{(k-1)}_{k-1}(\log X)\over (s-1)^k}\right),\cr}\eqno(3.4)
$$
which provides analytic continuation of the left-hand side of (3.4)
to $\CC$. We also have
$$
\int_X^\infty x^{-s}\d\D_k(x) = -X^{-s}\D_k(X) +
s\int_X^\infty x^{-s-1}\D_k(x)\d x.\eqno(3.5)
$$
Note that the last integral converges absolutely for $\s >\b_k$,
in view of the Cauchy-Schwarz inequality for integrals and
the definition (2.1) of $\b_k$. Therefore from (3.1)--(3.5) we obtain,
for $\max(\hf,\b_k) < \s \le 1$ and $T \le t \le 2T$,
$$
\eqalign{
\z^k(s) &= \sum_{n\le X}d_k(n)n^{-s} + s\int_X^\infty x^{-s-1}\D_k(x)\d x\cr&
+ O_\e\left(X^{\b_k-\s+\e} + T^{-1}X^{1-\s}\log^{k-1}X\right).\cr}\eqno(3.6)
$$
Observe now that (2.2) follows from
$$
\int_T^{2T}|\z(\s+it)|^{2k}\d t = \sum_{n=1}^\infty d_k^2(n)n^{-2\s}T
+ O_\e\left(T^{{2(1-\sigma)\over 1-\b_k}+\varepsilon}\right)\eqno(3.7)
$$
on replacing $T$ by $T2^{-j}\,(j \in \NN)$ and summing all the results.
To evaluate the integral in (3.7), we suppose that
$\max(\hf,\b_k) < \s < 1$, we   use (3.6) and
$$
\eqalign{
|a+b|^2 &\,= |a|^2 + |b|^2 + 2\R a{\bar b},\cr
a &:= \sum_{n\le {1\over2} X}d_k(n)n^{-s},\cr
b &:= \sum_{{1\over2} X < m\le X}d_k(m)m^{-s} \cr&
\,+ s\int_X^\infty x^{-s-1}\D_k(x)\d x
+ O_\e\left(X^{\b_k-\s+\e} + T^{-1}X^{1-\s}\log^{k-1}X\right).\cr}
$$
The reason of this splitting of the sum in two sums is to have
$m$ and $n$ differ by unity at least, which is expedient to have
in the integration that will follow. Now note that
we have, by the mean value theorem for Dirichlet polynomials
(see [4, Chapter 4]) and $d_k(n) \ll n^\e$,
$$\eqalign{
\int_T^{2T}|a|^{2}\d t &= T\sum_{n\le {1\over2} X}d_k^2(n)n^{-2\s}
+ O\Bigl(\sum_{n\le {1\over2} X}d_k^2(n)n^{1-2\s}\Bigr)\cr&
= T\sum_{n=1}^\infty d_k^2(n)n^{-2\s} + O_\e\left(
TX^{1-2\s+\e} + X^{2-2\s+\e}\right).\cr&}\eqno(3.8)
$$
To evaluate the mean square of $|b|$ we may proceed directly by
squaring out the modulus, or we may use Lemma 4 of [8], which says
that
$$
\int\limits_{T_1}^{T_2}\,\Bigl|\int\limits_\a^\b g(x)x^{-s}\d x\Bigr|^2\d t
\le 2\pi\int\limits_\a^\b g^2(x)x^{1-2\s}\d x
\quad(s = \s+it\,,T_1 < T_2,\,\a < \b)
$$
holds if $g(x)$ is a real-valued,
integrable function on $[\a,\b]$, a subinterval
of $[2,\,\infty)$, which is not necessarily finite. We shall obtain
$$
\eqalign{
&\int_T^{2T}|b|^{2}\d t  \ll_\e T^2\int_T^{2T}
\left|\int_X^\infty x^{-\s-it-1}\D_k(x)\d x\right|^{2}\d t \cr&
+ \,TX^{2\b_k-2\s+\e} + TX^{1-2\s+\e} + X^{2-2\s+\e}\cr&
\ll_\e T^2\int_X^\infty x^{-2\s-1}\D^2_k(x)\d x
+ TX^{2\b_k-2\s+\e} + TX^{1-2\s+\e} + X^{2-2\s+\e}\cr&
\ll_\e T^2X^{2\b_k-2\s+\e} + TX^{1-2\s+\e} + X^{2-2\s+\e}.\cr&}\eqno(3.9)
$$
We  have
$$
\eqalign{
&\int_T^{2T} a{\bar b}\d t  = \sum_{n\le{1\over2} X}d_k(n)n^{-\s}
\int_T^{2T}n^{-it}\Bigl\{\sum_{{1\over2} X < m\le  X}d_k(m)m^{-\s+it} \cr&
+ (\s-it)\int_X^\infty x^{-\s-1+it}\D_k(x)\d x + O_\e\left(X^{\b_k-\s+\e}
+ T^{-1}X^{1-\s+\e}\right)\Bigr\}\d t.
\cr}\eqno(3.10)
$$
By direct integration it is found that
$$\eqalign{&
\sum_{n\le{1\over2} X}d_k(n)n^{-\s}
\int_T^{2T}n^{-it}\sum_{{1\over2} X < m\le  X}d_k(m)m^{-\s+it}\d t \cr&
\ll \sum_{n\le{1\over2} X}d_k(n)n^{-\s}\sum_{{1\over2} X < m\le  X}
d_k(m)m^{-\s}\left|\log{m\over n}\right|^{-1}\cr&
\ll_\e X^{\e-\s}\sum_{n\le{1\over2} X}n^{-\s}\sum_{{1\over2} X < m\le  X}
\left({X\over m -n} + 1\right) \ll_\e X^{2-2\s+\e},
\cr}\eqno(3.11)
$$
on using the elementary inequality
$$
{1\over\log(1 + x)} \;\le\;1 + {1\over x}\qquad(x > 0).
$$
Similarly, by using the first derivative test (see [4, Lemma 2.1]),
we obtain
$$
\eqalign{
& \sum_{n\le{1\over2} X}d_k(n)n^{-\s}
\int_T^{2T}n^{-it}(\s-it)\int_X^\infty x^{-\s-1+it}\D_k(x)\d x \d t\cr&
\ll T\sum_{n\le{1\over2} X}d_k(n)n^{-\s}\int_X^\infty x^{-\s-1}|\D_k(x)|
\left|\log{x\over n}\right|^{-1}\d x \cr&
\ll_\e TX^{1+\b_k-2\s+\e},\cr}\eqno(3.12)
$$
where the interchange of the order of integration is justified by
absolute convergence. Therefore from (3.10)-(3.12) it follows that
$$
\int_T^{2T} a{\bar b}\d t \,\ll_\e\, T^\e\left(TX^{1+\b_k-2\s+\e}
 + X^{2-2\s}\right),\eqno(3.13)
$$
so that finally from (3.8)--(3.10) and (3.13) we obtain
$$\eqalign{&
\int_T^{2T}|\z(\s+it)|^{2k}\d t = \sum_{n=1}^\infty d_k^2(n)n^{-2\s}T\cr&
+ O_\e\left\{\left(T^\e(TX^{1+\b_k-2\s} + X^{2-2\s} + TX^{1-2\s}
+ T^2X^{2\b_k-2\s}\right)\right\}.\cr}\eqno(3.14)
$$
Now in (3.14) we set $X^{2-2\s} \asymp TX^{1+\b_k-2\s}$, namely
$$
X \= cT^{1\over1-\b_k},\eqno(3.15)
$$
where the constant $c > 0$ is chosen in such a way that (3.3) is satisfied.
With the choice (3.15) it is seen that (3.14) becomes (3.7),
and the proof of Theorem 1 is completed.
\medskip
{\bf Corollary 1.}
$$\eqalign{
\int_1^T|\z(\s+it)|^6\d t &= T\sum_{n=1}^\infty d_3^2(n)n^{-2\s}
+ O_\e\left(T^{3(1-\s)+\e}\right)\quad(\hf<\s<1),\cr
\int_1^T|\z(\s+it)|^8\d t &= T\sum_{n=1}^\infty d_4^2(n)n^{-2\s}
+ O_\e\left(T^{{16\over5}(1-\s)+\e}\right)\quad(\hf<\s<1),\cr
\int_1^T|\z(\s+it)|^{10}\d t &= T\sum_{n=1}^\infty d_5^2(n)n^{-2\s}
+ O_\e\left(T^{{40\over11}(1-\s)+\e}\right)
\quad(\hf<\s<1),\cr
\int_1^T|\z(\s+it)|^{12}\d t &= T\sum_{n=1}^\infty d_6^2(n)n^{-2\s}
+ O_\e\left(T^{4(1-\s)+\e}\right)\quad(\hf<\s<1).\cr&}
$$
The formulas follow from  Theorem 1 with the values $\b_3 = {1\over3},\,
\b_4 = {3\over8},\, \b_6 \le \hf$ (see [4]) and $\b_5 \le {9\over20}$
(see [20]).

\medskip
{\bf Remark 2.}
It transpires  that the Lindel\"of hypothesis ($\z(\s+it) \ll |t|^\e$
for $\s > \hf$) is equivalent to
$$
\int_1^T|\z(\s+it)|^{2k}\d t = \sum_{n=1}^\infty d_k^2(n)n^{-2\s}T
+ O_\e\left(T^{{4k(1-\s)\over k+1}+\e}\right)\;(\hf < \s \le 1,\, k \ge 2).
\eqno(3.16)
$$
Namely the Lindel\"of hypothesis implies (see [4, Chapter 13]) $\b_k =
(k-1)/(2k)$ for $k \ge 2$, in which case (3.16)
follows from (2.2) and Theorem 1.
Conversely, if (3.16) holds, then by [4, Lemma 7.1] we have,
for $T^\e \le H \le \hf T$ and $\hf < \s \le 1$,
$$
\eqalign{
|\z(\s+iT)|^{2k} &\ll 1 + \log T\int_{T-H}^{T+H}|\z(\s - {1\over\log T}
+ it)|^{2k}\d t\cr&
\ll_\e H\log T + T^{{4k(1-\s)\over k+1}+\e},\cr}
$$
which yields the Lindel\"of hypothesis on taking $H = \hf T$ and
letting $k\to\infty$.

\medskip

{\bf Remark 3.}
Other explicit results can be obtained  from Theorem 1 with the
bounds for $\b_k$ furnished by  [9], some of which are hitherto
the sharpest ones.
Our method of proof can be used to obtain a sharpening of the general
result proved in [12], since we did not use the Cauchy-Schwarz inequality
in estimating $\int_T^{2T}a{\bar b}\d t$, which was done in [12] and [13].
Namely we integrated directly the
expressions in question, which led to a sharper estimate than the one
that would have resulted from the application of the Cauchy-Schwarz
inequality.

\medskip
{\bf Proof of Theorem 2.} From the well-known Mellin inversion integral
$$
e^{-x} = {1\over2\pi i}\int_{c-i\infty}^{c+i\infty}
\G(w)x^{-w}\d w\qquad(c > 0,\,x > 0)
$$
we obtain
$$
\sum_{n=1}^\infty d_k(n)e^{-n/Y}n^{-s} = {1\over2\pi i}
\int_{2-i\infty}^{2+i\infty} Y^w\G(w)\z^k(s+w)\d w. \eqno(3.17)
$$
for $1 \ll Y \ll T^C\,(C>0)\,, T \le t \le 2T,\,\st < \s < 1$.
We move the line of integration
in (3.17) to $\R w = \st - \s$. In doing this we encounter the
pole $w = 1- s$ with residue $O(T^{-A})$ for any fixed $A>0$ in view
of Stirling's formula for the gamma-function. There is also the simple
pole at $w=0$ with residue $\z^k(s)$. Therefore from (3.17) it follows that
$$
\int_T^{2T}|\z(\s+it)|^{2k}\d t = \int_T^{2T}|F|^{2}\d t
+ \int_T^{2T}|G|^{2}\d t + 2\R \int_T^{2T}F{\bar G}\d t,
$$
where
$$\eqalign{
F &:= \sum_{n\le Y\log^2Y}d_k(n)e^{-n/Y}n^{-s},\cr
G &:= O\left(Y^{\st-\s}\int_{-\log^2T}^{\log^2T}
|\z(\st+it+iv)|^ke^{-|v|}\d v + T^{-A}\right).\cr}
$$
Consequently we have, by the mean value theorem for Dirichlet polynomials,
$$\eqalign{&
\int_T^{2T}|F|^{2}\d t = T\sum_{n\le Y\log^2Y}d_k^2(n)e^{-2n/Y}n^{-2\s}
+ O_\e(Y^{2-2\s+\e})\cr&
= T\sum_{n\le Y\log^2Y}d_k^2(n)n^{-2\s}
+ O_\e(TY^{-1}\sum_{n\le Y\log^2Y}d_k^2(n)n^{1-2\s}) + O(Y^{2-2\s+\e})\cr&
= T\sum_{n=1}^\infty d_k^2(n)n^{-2\s}  + O_\e(T^{1+\e}Y^{1-2\s} + Y^{2-2\s+\e}).
\cr}
$$
We also have, by the definition of $\st$,
$$\eqalign{
\int_T^{2T}|G|^2\d t &\ll Y^{2\st-2\s}\int_{-\log^2T}^{\log^2T}
\left(\int_T^{2T}|\z(\st+it+iv)|^{2k}\d t\right)\d v + T^{1-2A}\cr&
\ll_\e T^{1+\e}Y^{2\st-2\s},\cr}
$$
on taking $A$ sufficiently large. Finally by using the Cauchy-Schwarz
inequality we obtain
$$\eqalign{
\int_T^{2T}F{\bar G}\d t&\,\ll_\e\, T^\e(T + Y^{2-2\s})^{1/2}T^{1/2}Y^{\st-\s}
\cr&\,\ll_\e\, T^{1+\e}Y^{\st-\s} + T^{1/2+\e}Y^{1+\st-2\s}.\cr}
$$
Putting together all the estimates, replacing $T$ by $T2^{-j}\,(j \ge 1)$
and summing over $j$ we obtain
$$
R(k,\s;T) \;\ll_\e\; T^\e\left(TY^{1-2\s} + Y^{2-2\s} + TY^{\st-\s}
+ T^{1/2}Y^{1+\st-2\s}\right).
$$
Now we take
$$
Y \;=\; T^{1\over2-\st-\s}
$$
to obtain
$$
R(k,\s;T) \;\ll_\e\; T^\e\left(T^{2-2\s\over2-\st-\s} +
T^{4+\st-5\s\over2(2-\st-\s)}\right).
$$
On noting that the condition
$$
{4+\st-5\s\over2(2-\st-\s)} \;\le\; {2-2\s\over2-\st-\s}
$$
reduces to $\st \le \s$, which is certainly true, we obtain then
(2.3).

\medskip
{\bf Corollary 2.}
$$\eqalign{
\int_1^T|\z(\s+it)|^6\d t &= T\sum_{n=1}^\infty d_3^2(n)n^{-2\s}
+ O_\e\left(T^{{24(1-\s)\over17-12\s}+\e}\right)
\quad({\txt{7\over12}}<\s<1),\cr
\int_1^T|\z(\s+it)|^8\d t &= T\sum_{n=1}^\infty d_4^2(n)n^{-2\s}
+ O_\e\left(T^{{16(1-\s)\over11-8\s}+\e}\right)\quad({\txt{5\over8}}<\s<1),
\cr\int_1^T|\z(\s+it)|^{10}\d t &= T\sum_{n=1}^\infty d_5^2(n)n^{-2\s}
+ O_\e\left(T^{{120(1-\s)\over79-60\s}+\e}\right)
\quad({\txt{41\over60}}<\s<1),\cr
\int_1^T|\z(\s+it)|^{12}\d t &= T\sum_{n=1}^\infty d_6^2(n)n^{-2\s}
+ O_\e\left(T^{{14(1-\s)\over9-7\s}+\e}\right)\quad({\txt{5\over7}}<\s<1).
\cr&}$$
The above formulas follow  from
Theorem 2 with the values (see [4, Chapter 8])
$\s^*_3 \le {7\over12},\, \s^*_4 \le {5\over8},\, \s^*_5 \le {41\over60},\,
\s^*_6 \le {5\over7}$. They improve (1.5) and complement those
furnished by Corollary 1.

\medskip
\heading{\bf4. The mean value of the Rankin--Selberg series}
\endheading

\no
The arguments used in the proof of Theorem 1 and Theorem 2 are of a
general nature and can be adapted to obtain mean value results for a
wide class of Dirichlet series. Instead of working out the details
in the general case, which would entail various technicalities, we prefer
to conclude by considering two specific examples. In this section
we shall deal with the mean value of the
so-called Rankin--Selberg series (see R.A. Rankin  [15], [16])
$$
Z(s) := \z(2s)\sum_{n=1}^\infty|a(n)|^2n^{1-\k-s} = \sum_{n=1}^\infty
c_nn^{-s} \quad(\s > 1),
$$
and in Section 5 we shall consider the zeta-function attached to holomorphic
cusp forms. Here as usual $a(n)$ denotes the $n$-th Fourier coefficient
of a holomorphic cusp form $\f(z)$ of weight
$\k$ with respect to the full modular group $SL(2,\ZZ)$. We also suppose
that $\f(z)$ is a normalized eigenfunction for the Hecke operators
$T(n)$, so that $a(1) = 1$ and $a(n) \in \RR$. We have
(see [5], [11] and [13]) $c_n \ll_\e n^\e$,
$$
\sum_{n\le x}c_n^2 \,\ll_\e\, x(\log x)^{1+\e},\quad
\sum_{n\le x}c_n \;=\; Ax + \D(x,\f)\qquad(A > 0)
$$
with Rankin's classical estimate (see [15]) $\D(x,\f) \,\ll\, x^{3/5}$, and
$$
\int_1^X\D^2(x,\f)\d x \;\ll_\e\; X^{2+\e}.
$$
This means that analogously to (3.6) we have
$$
Z(s) = \sum_{n\le X}c_n + s\int_X^\infty \D(x,\f)x^{-s-1}\d x
+ O_\e(T^{-1}X^{1-\s} + X^{{1\over2}-\s+\e})\eqno(4.1)
$$
for $\hf < \s \le 1,\, T \le t \le 2T$, where $1 \ll X \ll T^C$ and
$X \,(\in [Y,2Y]$) satisfies (this is the analogue of (3.3))
$$
\D(X,\f) \;\ll_\e\;X^{{1\over2}+\e}.
$$
Then we write
$$
Z(s) :=  D + E
$$
with
$$
\eqalign{
D &:= \sum_{n\le{1\over2} X}c_nn^{-s},\cr E &:=
\sum_{{1\over2} X<m\le X}c_mm^{-s} + s\int_X^\infty \D(x,\f)x^{-s-1}\d x
+ O_\e(T^{-1}X^{1-\s} + X^{{1\over2}-\s+\e}),\cr}
$$
and consider
$$
\int_T^{2T}|Z(s)|^2\d t = \int_T^{2T}|D|^2\d t
+ \int_T^{2T}|E|^2\d t + 2\R\int_T^{2T}D{\bar E}\d t.
$$
Similarly as in the proof of Theorem 1 we find that
$$
\int_T^{2T}|D|^2\d t = T\sum_{n=1}^\infty c_n^2n^{-2\s}
+ O_\e\left\{(TX^{1-2\s} + X^{2-2\s})\log^{1+\e}X\right\},
$$
$$
\int_T^{2T}|E|^2\d t \,\ll_\e\, X^{2-2\s+\e} + T^2X^{1-2\s+\e},
$$
$$
\int_T^{2T}D{\bar E}\d t \,\ll_\e\,  X^{2-2\s+\e} + TX^{{3\over2}-\s+\e}.
$$
With the choice $X = bT^2$, where $b>0$ is a suitable constant, we obtain
$$
\int_T^{2T}|Z(\s+it)|^2\d t = T\sum_{n=1}^\infty c_n^2n^{-2\s}
+ O_\e(T^{4-4\s+\e}),
$$
which easily  gives then

\medskip
THEOREM 3. {\it For fixed $\s$ satisfying $\hf < \s\le 1$ we have}
$$
\int_1^{T}|Z(\s+it)|^2\d t = T\sum_{n=1}^\infty c_n^2n^{-2\s}
+ O_\e(T^{4-4\s+\e}).\eqno(4.2)
$$

\medskip
{\bf Remark 4.} The asymptotic formula (4.2) improves, for
${3\over4} < \s \le 1$, the result of K. Matsumoto [13] who proved
$$
\int_1^{T}|Z(\s+it)|^2\d t = T\sum_{n=1}^\infty c_n^2n^{-2\s}
+ R(\s,T)
$$
with
$$
R(\s,T) \,\ll_\e\, \Bigg\{\aligned T^{{5\over2}-2\s+\e}\qquad
&\left({\txt{3\over4}} < \s < {\txt{12+\sqrt{19}\over20}} =
0.81666\ldots\right),
\\T^{{60(1-\s)\over29-20\s}+\e} \qquad &\left({\txt{12+\sqrt{19}\over20}}
< \s < 1\right).\endaligned
$$
For $\hf < \s \le {3\over4}$ our result is slightly weaker
than the corresponding result of [13], namely
$$
R(\s, T) \,\ll_\e\, T^{4-4\s}(\log T)^{1+\e},
$$
but it should be remarked that (4.2) is a true asymptotic formula
only in the range ${3\over4} < \s \le 1$.
\bigskip
\medskip
\heading{\bf5. The mean value of the zeta-function of cusp forms}
\endheading
We retain the notation of Section 4 and consider (see [6]) the Dirichlet
series
$$
F(s) \;:=\; \sum_{n=1}^\infty {\tilde a}(n)n^{-s}\qquad(\s > 1),\eqno(5.1)
$$
which may be continued analytically to an entire function over $\CC$.
In (5.1) the arithmetic function
$$
{\tilde a}(n) \;:= a(n)n^{{1\over2}(1-\k)}\eqno(5.2)
$$
is the ``normalized" function of cusp form coefficients. This function
is ``small", since it satisfies $\ti \ll d(n)$ by Deligne's classical
estimate. We shall also consider
$$
F^2(s) \;=\; \sum_{n=1}^\infty \Ti n^{-s}\qquad(\s > 1),\eqno(5.3)
$$
where
$$
\Ti \;:=\; \sum_{d|n} {\tilde a}(d){\tilde a}({n\over d})
$$
is the  convolution of $\ti$ with itself. The mean values of $F(s)$ and
$F^2(s)$ were considered in [6]. It was proved there that, for $\s$ fixed,
$$
\int_1^T|F(\s + it)|^2\d t =T\sum_{n=1}^\infty |\ti|^2n^{-2\s}
+ H(\s;T)\eqno(5.4)
$$
and
$$
\int_1^T|F(\s + it)|^4\d t =T\sum_{n=1}^\infty |\Ti|^2n^{-2\s}
+ K(\s;T)\eqno(5.5)
$$
with
$$
H(\s;T) \;\ll_\e\; T^{{3\over2}-\s+\e} \qquad(\hf < \s<1)\eqno(5.6)
$$
and
$$
K(\s;T) \;\ll_\e\; T^{{11-8\s\over6}+\e} \qquad(\hf < \s<1).\eqno(5.7)
$$
Note that from (5.5) and (5.7) it transpires that we can obtain a true
asymptotic formula for the fourth moment of $F(\s+it)$ for ${5\over8}
< \s < 1$. This reflects the fact that we have (see [6])
$$
\int_1^T|F(\s + it)|^4\d t \;\ll_\e\; T^{1+\e}\eqno(5.8)
$$
only for $\s\ge {5\over8}$, and any improvement of the range for which
(5.8) holds would result in the  improvement of the  bound for $K(\s;T)$.
We shall improve on (5.6) and (5.7) by proving

\medskip
THEOREM 4. {\it If $H(\s;T)$ is defined by} (5.4), {\it then for $\s$ fixed
we have}
$$
H(\s;T) \ll_\e T^{{4(1-\s)\over3-2\s}+\e}\; (\hf < \s \le {\txt{3\over4}}),
\quad
H(\s;T) \ll_\e T^{{8\over3}(1-\s)+\e}\; ({\txt{3\over4}} \le \s \le 1).
\eqno(5.9)
$$

\medskip
THEOREM 5. {\it If $K(\s;T)$ is defined by} (5.5), {\it then for $\s$ fixed
we have}
$$
K(\s;T) \ll_\e T^{{16(1-\s)\over12-8\s}+\e}\; (\hf < \s \le {\txt{3\over4}}),
\quad
K(\s;T) \ll_\e T^{{16\over5}(1-\s)+\e}\; ({\txt{3\over4}} \le \s \le 1).
\eqno(5.10)
$$

\medskip
{\bf Proof of Theorem 4 and Theorem 5}. The first bounds in (5.9) and (5.10)
are the analogues of (2.3) of Theorem 2 corresponding to the values
$\s^*_1 = \hf$ and $\s_2^* = {5\over8}$, which follow from (5.4)
and (5.8), respectively. The method of proof of Theorem 2 may be used, since
$$
\ti \ll d(n) \ll_\e n^\e,\qquad \Ti \ll
\sum_{\delta|n}d(\delta)d({n\over\delta}) \ll_\e n^\e.\eqno(5.11)
$$
Similarly the second bounds in (5.9) and (5.10) are the analogues of (2.2)
of Theorem 1 corresponding to the values $\b_1 = {1\over4}$ and
$\b_2 = {3\over8}$, respectively. Namely if we define
$$
\rho = \inf\left\{\,c\ge0\;:\;\int_1^X\Bigl(\sum_{n\le x}\ti\Bigr)^2\d x
\ll X^{1+2c}\,\right\}
$$
and
$$
\t = \inf\left\{\,c\ge0\;:\;\int_1^X\Bigl(\sum_{n\le x}\Ti\Bigr)^2\d x
\ll X^{1+2c}\,\right\},
$$
then $\rho = {1\over4}$ (this corresponds to $\b_2 = {1\over4}$ in the
classical Dirichlet divisor problem) and $\t \le {3\over8}$  (see [6];
this corresponds to $\b_4 = {3\over8}$ in the
Dirichlet divisor problem for $\D_4(x)$). Thus proceeding as in the
proof of Theorem 1 and keeping in mind again that (5.11) holds, we shall
obtain the second bounds in (5.9) and (5.10). Clearly the bounds in
(5.9) and (5.10) improve (5.6) and (5.7), respectively. Although the
first bound in (5.10) holds for $\hf < \s \le {\txt{3\over4}}$, it
is relevant only in the range $\s > {5\over8}$, when
$16(1-\s)/(11-8\s) < 1$, when (5.5) becomes a true asymptotic formula.

\bigskip
\vfill
\break
\topglue1cm

\Refs

\bigskip

\item{[1]} R. Balasubramanian, A. Ivi\'c and K. Ramachandra,  An
application of the Hooley--Huxley contour, {\it Acta Arith}. {\bf65}(1993),
45-51.

\item{[2]} K. Chandrasekharan and R. Narasimhan,  Functional
equations with multiple gamma factors  and the average order of
arithmetical functions, {\it Annals Math.} {\bf76}(1962), 93-136.

\item{[3]} K. Chandrasekharan and R. Narasimhan,
The approximate functional equation for a class of
zeta-functions, {\it Math. Annalen} {\bf 152}(1963), 30-64.

\item{[4]} A. Ivi\'c,  The Riemann zeta-function, {\it John
Wiley \& Sons}, New York, 1985.

\item{[5]} A. Ivi\'c,  Large values of certain number-theoretic
error terms, {\it Acta Arith.} {\bf56}(1990), 135-159.

\item{[6]} A. Ivi\'c,  On zeta-functions asociated with
Fourier coefficients of cusp forms, {\it Proceedings of the Amalfi Conference
on Analytic Number Theory (Amalfi, September 1989)}, Universit\`a di Salerno,
Salerno 1992, 231-246.

\item{[7]} A. Ivi\'c,  Some problems on mean values
of the Riemann zeta-function, {\it J. de Th\'eorie des Nombres de Bordeaux}
{\bf8}(1996), 101-123.

\item{[8]} A. Ivi\'c, On some conjectures and results for the
Riemann zeta-function and Hecke series, {\it Acta Arithmetica}
{\bf99}(2001), 115-145.

\item{[9]} A. Ivi\'c and M. Ouellet, Some new estimates in the
Dirichlet divisor problem, {\it Acta Arith.} {\bf52}(1989), 241-253.

\item{[10]} A. Ivi\'c and K. Matsumoto,
On the error term in the
mean square formula for the Riemann zeta-function in the critical strip,
{\it Monatshefte Math.} {\bf121}(1996), 213-229.

\item{[11]} A. Ivi\'c, K. Matsumoto and Y. Tanigawa,  On Riesz
means of the coefficients of the Rankin--Selberg series,
{\it Math. Proc. Camb. Phil. Soc.} {\bf127}(1999), 117-131.

\item{[12]} S. Kanemitsu, A. Sankaranarayanan and Y. Tanigawa,
A mean value theorm for Dirichlet series and a general divisor
problem, Monats. Math. {\bf136}(2002), 17-34.

\item{[13]} K. Matsumoto,  The mean values and the universality
of Rankin-Selberg L-functions, {\it Proc. Turku Symposium on
Number Theory in Memory of K. Inkeri}, May 31--June 4, 1999,
Walter de Gruyter, Berlin etc., 2001, pp. 201-221.

\item{[14]} A. Perelli,  General L-functions, {\it Ann. Mat. Pura
Appl.} {\bf130}(1982), 287-306.

\item{[15]} R.A. Rankin,  Contributions to the theory of Ramanujan's
function $\tau(n)$ and similar arithmetical functions. II, The order
of Fourier coefficients of integral modular forms, {\it Math. Proc.
Cambridge  Phil. Soc.} {\bf 35}(1939), 357-372.

\item{[16]} R.A. Rankin,  Modular Forms, {\it Ellis Horwood Ltd.,}
Chichester, England, 1984.

\item{[17]} H.-E. Richert,  \"Uber Dirichletreihen mit
Funktionalgleichung, {\it Publs. Inst. Math. (Belgrade)}
{\bf11}(1957), 73-124.

\item{[18]} A. Selberg,   Old and new conjectures and results
about a class of Dirichlet series, {\it Proc. Amalfi Conf. Analytic
Number Theory, eds. E. Bombieri et al.,} Universit\`a di Salerno,
Salerno, 1992, 367-385.

\item{[19]} E.C. Titchmarsh,  The theory of the Riemann
zeta-function, 2nd edition,{ \it Oxford University Press}, Oxford, 1986.

\item{[20]} Zhang Wenpeng, On the divisor problem, {\it Kexue Tongbao}
{\bf33}(1988), 1484-1485.

\vskip2cm

Aleksandar Ivi\'c

Katedra Matematike RGF-a

Universitet u Beogradu

\DJ u\v sina 7, 11000 Beograd

Serbia and Montenegro

{\tt aivic\@rgf.bg.ac.yu,\enskip aivic\@matf.bg.ac.yu}

\endRefs

\bye